%% file: ifacconf.tex
\documentclass{ifacconf}

\usepackage{graphicx}      
\usepackage{natbib}        

\usepackage{xcolor}
\usepackage{amsmath}
\usepackage{amsfonts} 

\usepackage{multirow}
\usepackage{booktabs}

\usepackage{import}
\begin{document}
\begin{frontmatter}


\title{Data-driven observer design for an inertia wheel pendulum with static friction\thanksref{footnoteinfo}}

\thanks[footnoteinfo]{This work has been supported by the COMET-K2 Center of the Linz Center of Mechatronics (LCM) funded by the Austrian federal government and the federal state of Upper Austria.
}

\author[First]{L. Ecker}
\author[First]{M. Schöberl} 

\address[First]{Institute of Automatic Control and Control Systems Technology, Johannes Kepler University Linz,  Altenberger Str. 69, A-4040 Linz, Austria (e-mail: \{lukas.ecker, markus.schoeberl\}@jku.at).}

\begin{abstract} 
	An indirect data-driven state observer design approach for the inertia wheel pendulum considering static friction of the actuated inertia disc is presented. 
	The frictional forces occurring in a real laboratory setup are characterized by the Stribeck effect as well as the transition between two different dynamic behaviors, sticking and non-sticking. These switching nonlinear dynamics are identified with various machine learning methodologies in a data-driven manner, i.e., the unsupervised separation and feature clustering of measured state trajectories into two dynamic classes, and  the supervised classification of a state-dependent switching condition. The identified system with the interior switching-structure of two dynamics is combined with a moving horizon estimator.
\end{abstract}

\begin{keyword}
	Switched systems, System identification, Machine learning, Statistical data analysis, Data-based observer design, Inertia wheel pendulum, Static friction
\end{keyword}

\end{frontmatter}

\section{Introduction}

The mathematical modeling of a considered dynamic process is usually the first step of any control engineering task, whether it is a system theoretical investigation or a controller design. The modeling is in general a very time-consuming process, whose effort increases dramatically with the complexity and the desired accuracy, see \cite{Nelles:2001}. The complexity, such as occurring nonlinearities or structural switchovers, cannot be influenced. However, a modest model accuracy is often sufficient for various applications, although, this is in general not the case for the observer design.
The validity of the model is here crucial for the design and performance, since deviations between the model and the real process can hardly be compensated by the observer, causing biased state estimates.
{In addition to the derivation of the governing differential equations in the modeling process, unknown model parameters have to be determined experimentally.}
Data-driven approaches intend to minimize the modeling effort while using only recorded data. This is either done by directly determining the pursued objective, i.e., a control law or an observer, or by first generating a model and then using it to indirectly solve for the objective, see, e.g., \cite{Allgoewer:2021} and \cite{Brunton:2019}.

The present manuscript addresses an exemplary indirect data-driven observer design approach for a nonlinear model with structure switching. The problem statement is discussed in Section 2. Section 3 is dedicated to the data-driven modeling step while in Section 4 the actual observer design process is explained and experimental results are presented.



\begin{figure}[h]
	\begin{center}
		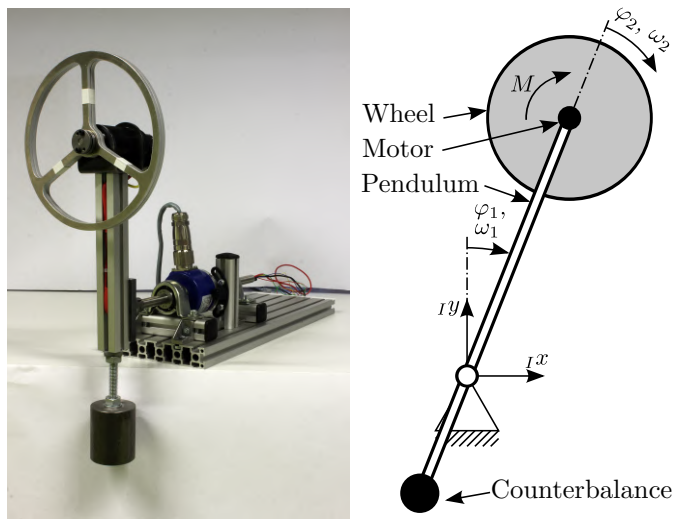
		\caption{Inertia wheel pendulum}
		\label{fig:inertia_wheel_pendulum}
	\end{center}
\end{figure}


\section{Problem statement}
	{
	The focus of this paper is on the data-driven system identification and observer design for the nonlinear switched model of the inertia wheel pendulum with static friction. The inertia wheel pendulum, as depicted in Fig.  \ref{fig:inertia_wheel_pendulum}, is a common control benchmark problem, see, e.g., \cite{Spong:1999}. In order to compare the subsequent results of the proposed data-driven approaches, the mathematical model of the pendulum is first briefly discussed. A more detailed mathematical investigation of the considered pendulum with static friction and its switching condition can be found in \cite{Ecker:2022b}, where, in contrast to the data-driven approaches in this manuscript, three observers were designed on the basis of the mathematical model.}
	
	
\subsection{Inertia wheel pendulum}
	
	The inertia wheel pendulum is composed of a rigid, ball-bearing pendulum, a motor-driven wheel and a counterweight. The absolute angle and its corresponding angular velocity of the pendulum are specified by the coordinates $\varphi_{1}$ and $\omega_{1}$, respectively. The inertia wheel, described by the relative angle $\varphi_{2}$ and the angular velocity $\omega_{2}$, is driven by a subordinate torque-controlled DC electric motor. Since the electrical dynamics are negligible compared to the mechanical dynamics, the motor torque $M$ is considered as the input to the system.
	The coordinate $\varphi_{2}$ will be omitted in the following considerations, since the actual angle of the wheel is not of interest for further control tasks.	
    The equations of motion for this simplified model can be derived by determining kinetic and potential energies and {applying} the Euler-Lagrange equations. However, as measurements on a real laboratory demonstrator and the validations in \cite{Ecker:2022b} show, the behavior of the system is affected by two dynamics and non trivial frictional forces. The two dynamic classes of the switched system are abbreviated as $\mathcal{C}_1$ and $\mathcal{C}_2$, which in turn describe the pendulum with sticking and with non-sticking wheel, respectively. An associated switching condition $g(x,u)$ partitions the state space $ x = [\varphi_{1}, \omega_{1}, \omega_{2}]^T $ and input space $ u = M$ with $(x,u) \in \mathcal{R}_1 \sqcup \mathcal{R}_2 $ into disjoint decision regions $\mathcal{R}_{1}$ and $\mathcal{R}_{2}$, such that all elements in $\mathcal{R}_j$ with $j \in \{1,2\}$ are assigned to class $\mathcal{C}_j$. In the considered scenario, the state-dependent sticking condition, as analytically derived in \cite{Ecker:2022b}, is given by
    \makeatletter
    \renewcommand{\maketag@@@}[1]{\hbox{\m@th\normalsize\normalfont#1}}%
    \makeatother
	\begin{equation}\small
		g(x,u) = \left\{ \begin{array}{ll}
		\mathcal{C}_1   & \left|\frac{\theta_{2}}{\theta_{1}+\theta_{2}}\left( d_1 \omega_1-a \sin(\varphi_{1})\right)+ M \right| < r_S \\
		\mathcal{C}_2 &\text{otherwise}
		\end{array} \right. . \label{equ:switching_condition_ref}
	\end{equation}
	A summary of the occurring model quantities is given in Table \ref{tab:parameter}.
	Depending on the current state and its class affiliation, the evolution of the state $\dot{x} = f_{\mathcal{C}_j}(x,u)$ is governed by a sticking dynamic $f_{\mathcal{C}_1}=[
	\omega_{1},	\frac{1}{\theta_{1}} \left( a\sin \left( \varphi_{1} \right) -d_{1}\omega_{1} \right),	0]^T$
	or	a non-sticking dynamic
	\begin{equation*}
	f_{\mathcal{C}_2}
	=
	\begin{bmatrix}
	\omega_{1}\\
	\frac{1}{\theta_{1}} \left( a\sin \left( \varphi_{1} \right) -d_{1}\omega_{1} +d_{2}\omega_{2}-M {+ M_S}\right)\\
	-\frac{a}{\theta_{1}}\sin(\varphi_{1}) + \frac{d_1}{\theta_{1}}  \omega_{1} +\frac{1}{\theta_{c}} (M {- M_S} - d_2 \omega_{2})
	\end{bmatrix}.
	\end{equation*}
	
	In addition to the viscous damping terms of the motor shaft and the ball bearings of the pendulum, which already appear in the equations of motion as $d_1 \omega_{1}$ and $d_2 \omega_{2}$, the static friction of the wheel, denoted $M_S$, was modeled according to the Stribeck effect as 
	\begin{equation}\label{eq:sticking_condition}
		M_S = r_C \mathrm{sign}(\omega_2) + (r_S - r_C) \exp^{- \left(\frac{\omega_2}{\omega_{2,0}}\right)} \mathrm{sign}(\omega_{2})
	\end{equation}
	
	\begin{table}[ht]
		\begin{center}
			
			\caption{Parameters of the inertia wheel pendulum and the static friction model.}\label{tab:parameter}
			{
				\begin{tabular}{llll}
					\toprule
					
					\multicolumn{2}{l}{Model parameters} & \multicolumn{2}{l}{Static friction model} \\ 
					\cmidrule{1-4}
					\multicolumn{1}{l}{Sym.} &  Description  &  \multicolumn{1}{l}{Sym.}  & Description \\
					\midrule 
					
					$a$ & Potential terms &  $r_C$& Coulomb friction \\
					$\theta_{1}$ &  Moment of
					inertia pendulum & $r_S$ &Static friction\\
					$\theta_{2}$ &  Moment of inertia wheel
					 & $w_{2,0}$&Stribeck velocity \\
					$d_1$ & Viscous damping pendulum&  \\
					$d_2$ & Viscous damping wheel&  \\
					\bottomrule
			\end{tabular}}
		\end{center}
		
	\end{table}

\subsection{System identification and observer design}

	As briefly motivated in the introduction section, the main objective is to propose an indirect data-driven observer design, i.e., the entire mathematical modeling should be replaced by an accurate system identification capable of handling different dynamics with a state-dependent switching condition. Subsequently, the state observer is designed on the basis of the identified model with the pendulum angle $\varphi_{1}$ as measurable output.
	The main problems that arise here, in addition to the fundamental difficulties with the nonlinear observer design of a switched system, see, e.g., \cite{Liberzon:2003} and \cite{Simon:06}, are mainly due to the intended purely data-driven modeling. {With no prior knowledge and only recorded state and input trajectories, the dynamics should be identified and the switching condition approximated.}
	Both the identification step and the observer design are based on sampled-data systems, where $k$ corresponds to the time index of the discrete-time state $x_k$ and $T_a$ corresponds to the sampling time of the recorded data sets. Note that also an identification of a continuous-time model would be conceivable.  While this would allow the use of different quasi-continuous observer design {methodologies}, it has a slightly negative aspect for the accuracy of the model obtained in the identification step, since additional numerical derivatives of recorded measurements are required, see, e.g., \cite{Chartrand:2011}.
	
\begin{figure}[h]
		\begin{center}
			\def\svgwidth{8.4cm}
			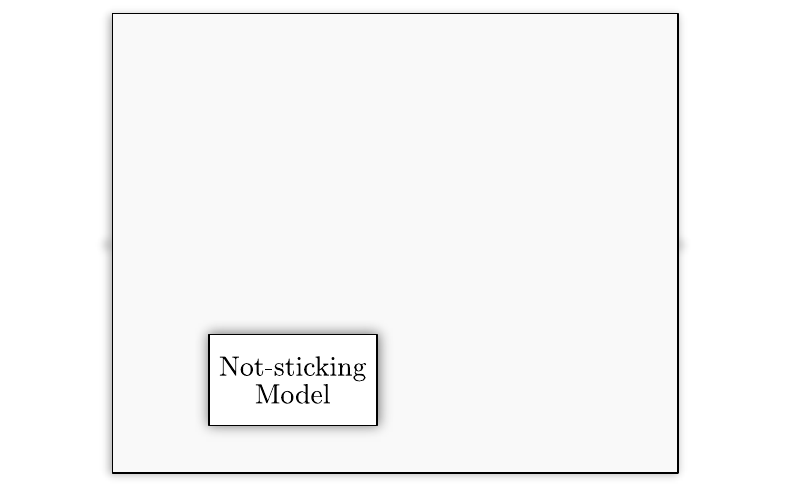

			\caption{Schematic of the inertia wheel pendulum with switching structure.} 
			\label{fig:iwp_switching_structure}
		\end{center}
	\end{figure}

	However, the identification step can be divided into several subtasks using different machine learning methods.
	In a first step, the task is to cluster the individual time series tuples consisting of state $x_i$, input $u_i$ and successor state $x_{i+1}$ for $i = 0 \ldots N-1$, where $N$ corresponds to the number of recorded samples, into two classes $\mathcal{C}_1$ and $\mathcal{C}_2$. Each class represents a different dynamic of the inertia wheel pendulum. The purpose of this unsupervised assignment is to generate target labels that subsequently allow the separate identification of the sticking and the non-sticking model, as well as the classification of the switching condition. In order to determine the nonlinear governing equations of the inertia wheel pendulum, spare regression methods are used for the identification of two parsimonious models, see \cite{Brunton:2016}. 
	{As mentioned, the clustering step could be omitted and only one model could be identified, but this would drastically affect the quality of the overall model and thus also the performance of the observer.}
	Moreover, the final composition of the overall model requires a clear definition of the decision regions $\mathcal{R}_1$ and $\mathcal{R}_2$.
	To this end, the previously clustered target labels are used in a supervised machine-learning environment to approximate the corresponding state-dependent sticking condition $\hat{g}(x_k, u_k)$. The principle schematic of the identified inertia wheel pendulum is illustrated in Fig. \ref{fig:iwp_switching_structure}.
	
\section{Data-driven system modeling}

This section is dedicated to the data-driven system modeling of the inertia wheel pendulum with static friction. For this purpose, this part is organized into three subsections covering the three main mentioned steps of the modeling process. First, the generation and unsupervised clustering of the training data is presented.
In the second stage, a classifier is trained with the previously determined clusters to approximate the state-dependent switching condition.
Finally, two models, one for the sticking dynamic and one for the non-sticking dynamic, are identified and linked via the approximated switching condition from step two.

\subsection{Dynamical feature clustering}

\begin{figure}[h]
	\begin{center}
		\includegraphics[width=8.4cm]{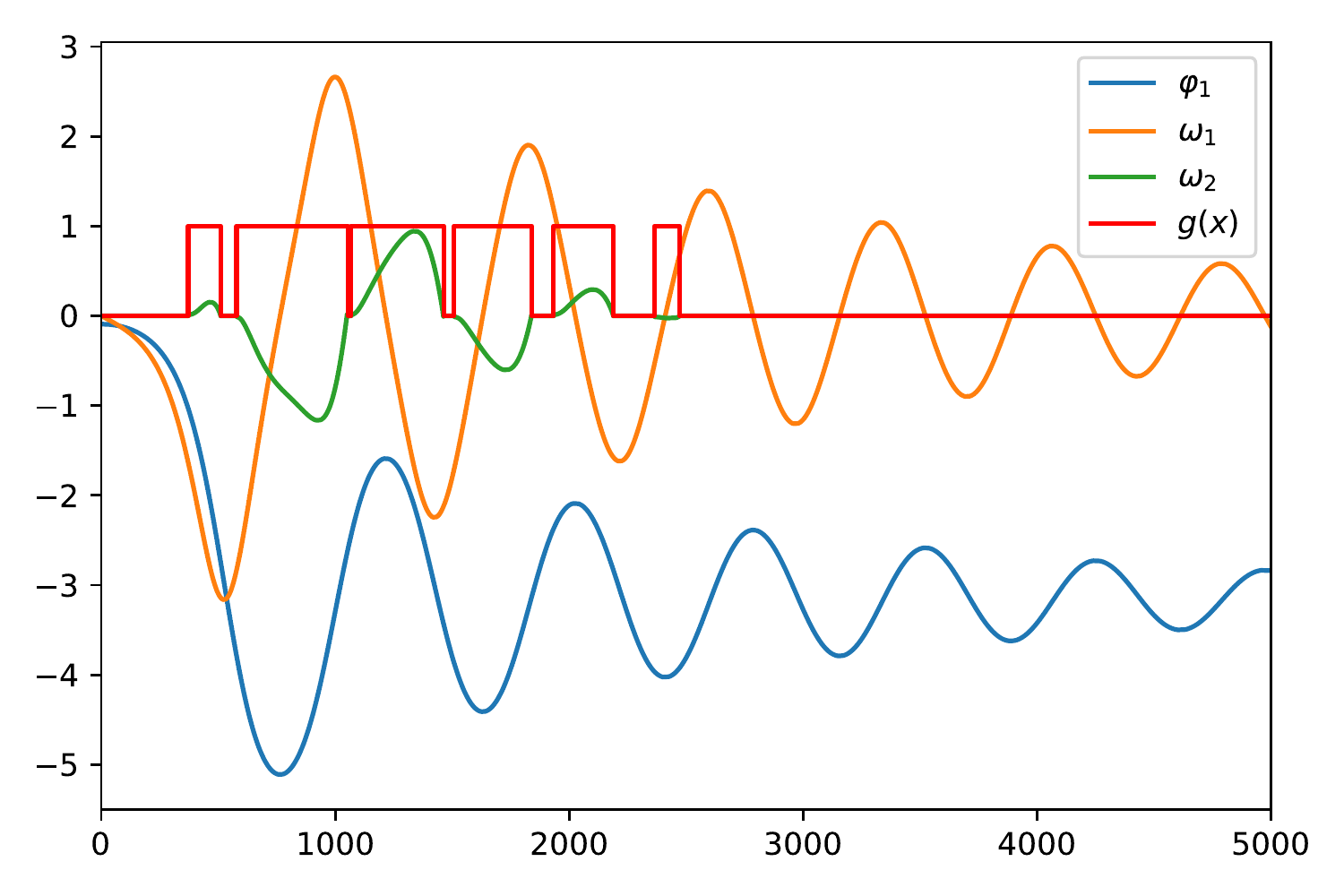}
		\caption{Drop-down experiment with evaluated sticking-condition reference.} 
		\label{fig:iwp_drop_down_experiement}
	\end{center}
\end{figure}

The training data set for the overall identification process consists of approximately 125000 recorded state and input samples from unactuated pendulum drop down experiments and recordings with random torque step sequences. Since the focus is on the behavior of the data-driven approach in combination with the switching system, the results of dropdown experiments, as exemplified in Fig. \ref{fig:iwp_drop_down_experiement}, are examined. Note that for the purpose of reference, the evaluation of the sticking-condition $g(x_k, u_k)$, cf. Eq. (\ref{equ:switching_condition_ref}), is plotted in addition to the state trajectories with value zero corresponding to class $\mathcal{C}_1$ and one to $\mathcal{C}_2$. This information is of course not available for the data-driven modeling.
The training sets with a non-zero input are obviously necessary for the identification of the input behavior. In contrast, the drop down experiments with zero input are necessary for the identification of the switching condition and the sticking model, since already small input torques release the disc from adhesion. The number of recordings from both types of experiments should be balanced in order to avoid artificially introduced biases that could affect the used probabilistic methods. The main purpose of this step is to cluster the set of data tuples $\{(x_i,u_i,x_{i+1})\}_{i=0\ldots N-1}$, i.e., to assign each sample tuple $(x_i,u_i,x_{i+1})$ to either class $C_1$ or $C_2$. The assignment is intended to separate the samples according to their dynamic transition from $x_k$ and $u_k$ to $x_{k+1}$. Therefore, it is necessary that both dynamics are geometrically separable in the feature space of $(x_i,u_i,x_{i+1})$.
In general, there are several algorithms from the field of unsupervised machine learning, such as k-means, spectral clustering, OPTICS, DBSCAN, and various mixture model approaches that can be used for this task, see, e.g., \cite{Bishop:2006}. 

\begin{figure}[h]
	\begin{center}
		\includegraphics[width=8.4cm]{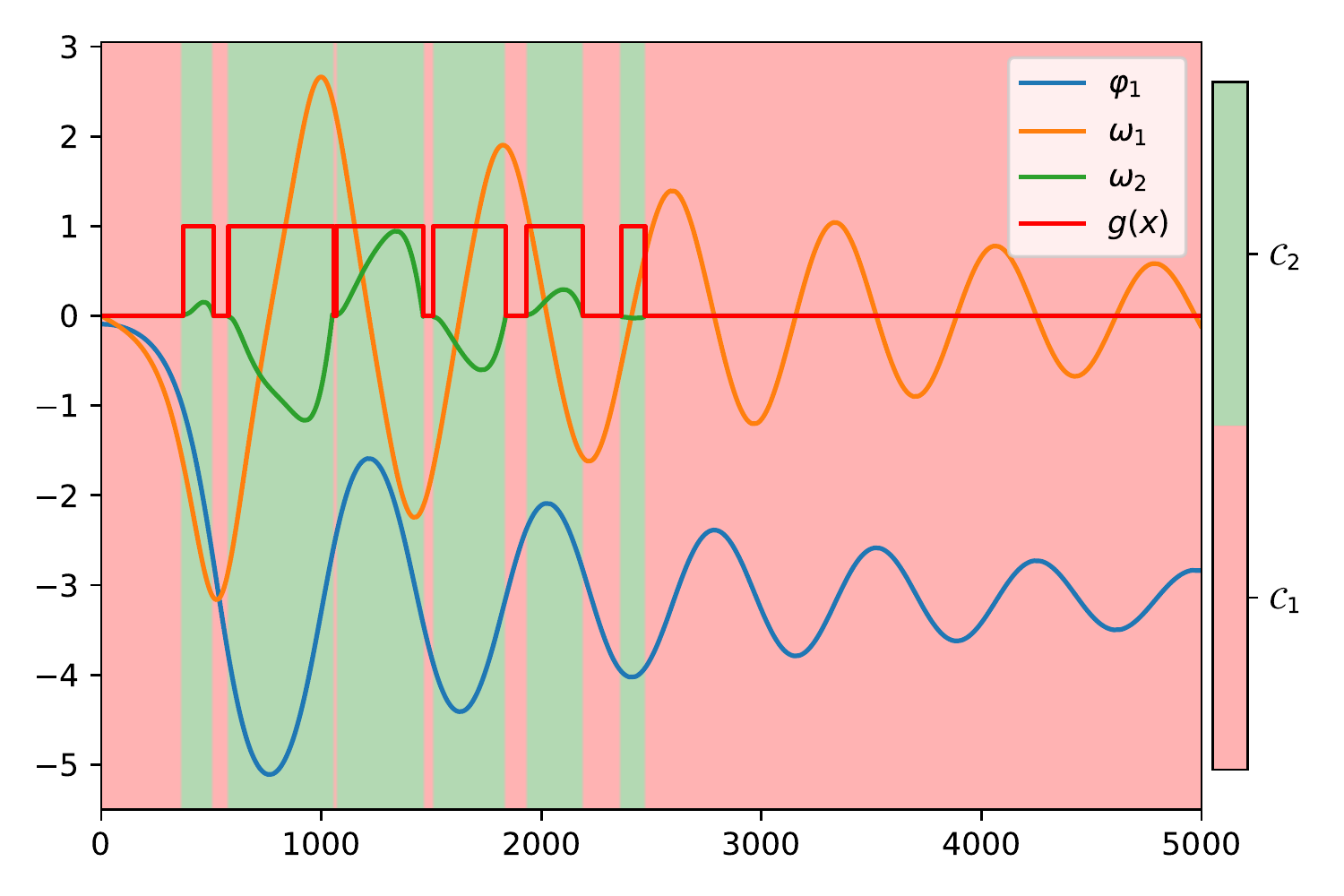}
		\caption{Drop-down experiment with highlighted clusters.} 
		\label{fig:iwp_drop_down_experiement_clustered}
	\end{center}
\end{figure}

However, in this manuscript, a two-component Gaussian mixture model is trained with an expectation maximization algorithm, see \cite{Gupta:2011}. This probabilistic model assumes that the samples of the classes $\mathcal{C}_1$ and $\mathcal{C}_2$ are distributed according to two multivariate Gaussian probability densities. The corresponding unknown first two stochastic moments can then be determined as maximum likelihood estimates using expectation maximization. The probability density function of the two-component Gaussian mixture distribution for the sample $\xi_i = (x_i,u_i,x_{i+1})$ can be written as $$p(\xi_i) = \pi_{\mathcal{C}_1} \mathcal{N}(\xi_i|\mu_{\mathcal{C}_1},\Sigma_{\mathcal{C}_1}) + \pi_{\mathcal{C}_2} \mathcal{N}(\xi_i|\mu_{\mathcal{C}_2},\Sigma_{\mathcal{C}_2})$$
with  $\pi = [\pi_{\mathcal{C}_1}, \pi_{\mathcal{C}_2}]$ denoting the mixing coefficients and $\mu = [\mu_{\mathcal{C}_1},\mu_{\mathcal{C}_2}]$ and $\Sigma = [\Sigma_{\mathcal{C}_1} > 0 , \Sigma_{\mathcal{C}_2} > 0]$ denoting the stochastic moments mean and variance of the multivariate Gaussian probability density functions $\mathcal{N}(\mu_{\mathcal{C}_j},\Sigma_{\mathcal{C}_j})$. 
In order to be a valid probability density function, $\pi_{\mathcal{C}_1} + \pi_{\mathcal{C}_2} = 1$ with $0\le \pi_{\mathcal{C}_j} \le 1$ must hold, see \cite{Papoulis:2002}. The mentioned unknown quantities are determined as the maximizing parameters of the log-likelihood function
\begin{equation}
\ln p(\{\xi_i\}_{i=0}^{N-1}|\pi,\mu,\Sigma) = \sum_{i=0}^{N-1} \ln \sum_{j=1}^{2} \pi_{\mathcal{C}_j} \mathcal{N}(\xi_i|\mu_{\mathcal{C}_j},\Sigma_{\mathcal{C}_j}).
\label{eq:log_likelihood_function}
\end{equation}
Due to the summation of the Gaussian mixtures, the derivatives of the log-likelihood function set to zero do not yield closed-form solutions. The derivative of Eq. (\ref{eq:log_likelihood_function}) with respect to the mean $\mu_{\mathcal{C}_j}$ yields
\begin{equation}
\mu_{\mathcal{C}_j} = \frac{\sum_{i=0}^{N-1}  p(\mathcal{C}_j|\xi_i) \xi_i}{\sum_{i=0}^{N-1} p(\mathcal{C}_j|\xi_i)}, \label{equ:em_mu}
\end{equation}
where $p(\mathcal{C}_j|\xi_i) = \frac{\pi_j \mathcal{N}(\xi_i|\mu_{\mathcal{C}_j},\Sigma_{\mathcal{C}_j})}{\sum_{l=1}^{2}\pi_l \mathcal{N}(\xi_i|\mu_{\mathcal{C}_l},\Sigma_{\mathcal{C}_l})}$  corresponds to the a posteriori probability of class $\mathcal{C}_j$ given $\xi_i$. Likewise, the derivative of Eq. (\ref{eq:log_likelihood_function}) with respect to the variance $\Sigma_{\mathcal{C}_j}$ results in 
\begin{equation}
\Sigma_{\mathcal{C}_j} = \frac{ \sum_{i=0}^{N-1} p(\mathcal{C}_j|\xi_i) (\xi_i - \mu_{\mathcal{C}_j})(\xi_i - \mu_{\mathcal{C}_j})^T}{\sum_{i=0}^{N-1} p(\mathcal{C}_j|\xi_i)}.  \label{equ:em_sigma}
\end{equation}
In order to take the constraint $\pi_{\mathcal{C}_1}+\pi_{\mathcal{C}_2} = 1$ properly in to account, we intend to maximize the function $\ln p(\{\xi_i\}_{i=0}^{N-1}|\pi,\mu,\Sigma) + \lambda \left(\pi_{\mathcal{C}_1} +  \pi_{\mathcal{C}_2} - 1\right)$, which leads after elimination of the Lagrange multiplier $\lambda$ to
\begin{equation}
\pi_{\mathcal{C}_j} = \frac{1}{N} \sum_{i=0}^{N-1} p(\mathcal{C}_j|\xi_i).  \label{equ:em_pi}
\end{equation}
The equations (\ref{equ:em_mu}), (\ref{equ:em_sigma}), and (\ref{equ:em_pi}) do not represent a closed solution, since $p(\mathcal{C}_j|\xi_i)$ in turn depends on the parameters $\mu$, $\Sigma$ and $\pi$. However, an iterative evaluation of these equations with initially defined parameter values complies to the special case of the expectation maximization algorithm for Gaussian mixture models. The initial values can be estimated, for instance, with a k-means algorithm and then used in the first expectation maximization step to calculate the posteriori probabilities $p(\mathcal{C}_j|\xi_i)$. These probabilities are then used in the second step to update the parameters of the mixture model according to the equations (\ref{equ:em_mu}), (\ref{equ:em_sigma}), and (\ref{equ:em_pi}). The iterative process of these two steps guarantees that the likelihood function (\ref{eq:log_likelihood_function}) increases in each iteration, see, e.g., \cite{Bishop:2006} and \cite{Gupta:2011}.
After the training of this probabilistic model, each sample $\xi_i$ can be assigned to its most likely class affiliation by evaluating $p(C_j|\xi_i)$ for $\mathcal{C}_1$ and $\mathcal{C}_2$. The obtained target label is thus determined as $$\gamma_i = \underset{C_j \in \{\mathcal{C}_1,\mathcal{C}_2\}}{\arg \max} p(C_j|\xi_i)$$ for each recorded feature sample $\xi_i$.
{The results of the clustering step are presented in Fig. \ref{fig:iwp_drop_down_experiement_clustered}. The particular feature clusters are illustrated with different background colors, i.e., the red background highlights class $\mathcal{C}_1$ and the green background highlights class $\mathcal{C}_2$.}
The transitions of the clusters coincide almost exactly with the plotted reference switching condition, cf. Eq (\ref{equ:switching_condition_ref}). This high consistency is very important for both the identification and classification steps, since otherwise the misclassified labels would distort the training. However, the results show that the mixture model is quite capable of clustering the data, although the assumption of a Gaussian distribution of the samples is not satisfied due to the inherent nonlinearities of the system under consideration. Further research will address this issue in more detail.

\subsection{Classification of the switching condition}
\begin{figure}[h]
	\begin{center}
		\includegraphics[width=8.4cm]{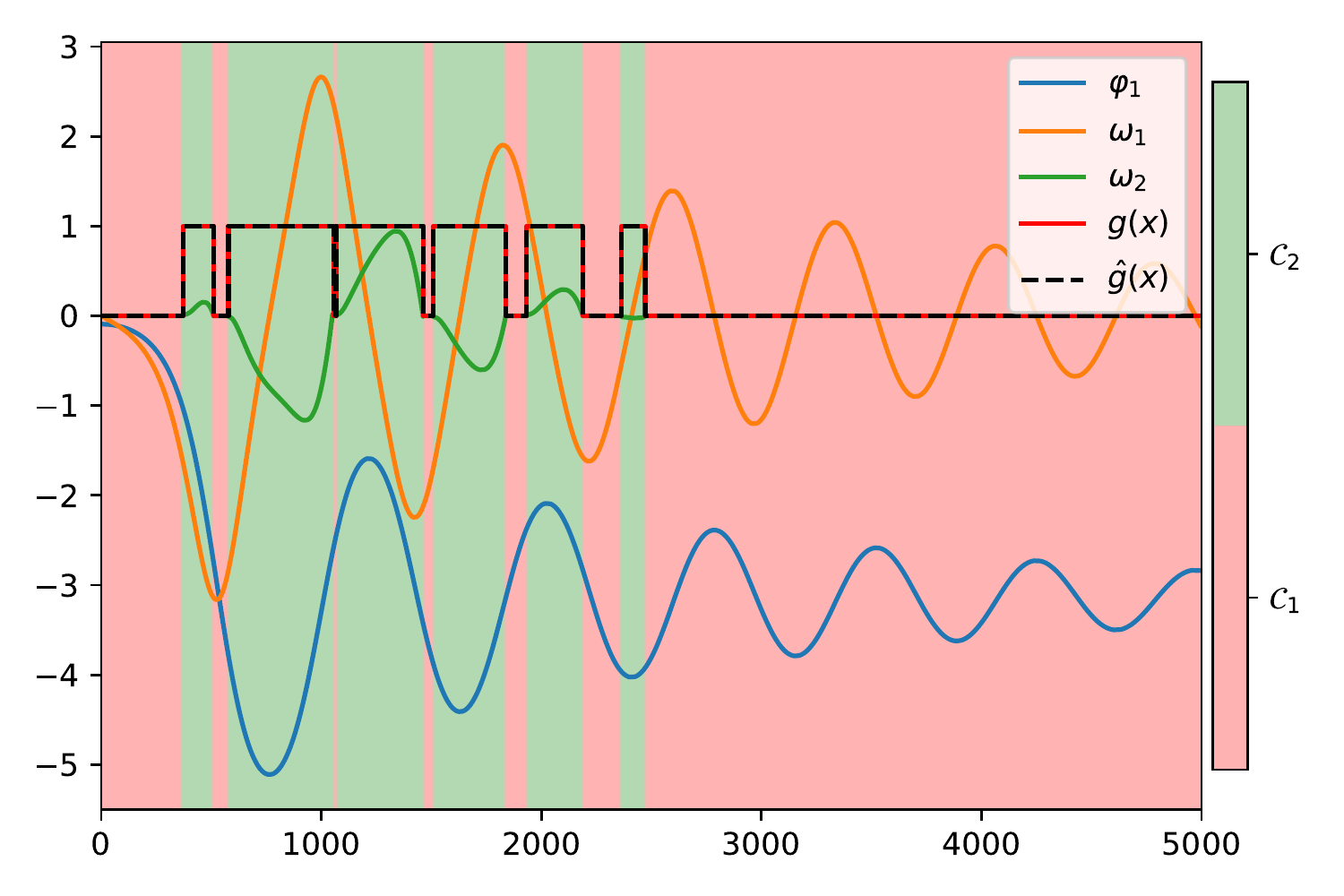}
		\caption{Drop-down experiment with approximated switching condition. }
		\label{fig:iwp_drop_down_experiement_classified}
	\end{center}
\end{figure}

The task in the previous section was to separate the training data sets into two classes. In this way, each sample was assigned a label, either $\mathcal{C}_1$ or $\mathcal{C}_2$. The purpose of this section is now to exploit this obtained information for the approximation of the switching condition and thus the partition of the state space in $\mathcal{R}_1$ and $\mathcal{R}_2$. 
In machine learning, this corresponds in principle to a supervised learning task where various classification algorithms such as logistic or naive Bayesian regression methods can be considered, see \cite{Bishop:2006}. In this manuscript, a decision tree learning algorithm is proposed for the classification of the switching condition, since decision trees are simple to interpret and have low computational costs for predicting data, see \cite{Breiman:1984}. In contrast to the clustering step, the considered feature space with training data tuples $\xi_i = (x_i,u_i)$ do not contain the successor state, since an approximation of the switching condition $g(x_k, u_k)$ with arguments $x_k$ and $u_k$ is intended. However, the target labels $\gamma_i  \in \{\mathcal{C}_1, \mathcal{C}_2\}$ from the previous clustering step are available.
A decision tree consists in general of multiple nodes and leaves. In this context, nodes represent intersections at which the feature space is partitioned with respect to a single feature in order to group feature samples with similar labels. Each node can lead either to a new node or to a leaf.  The end of the tree is represented by leaves on which the remaining feature space is assigned to class $\mathcal{C}_1$ or $\mathcal{C}_2$. Tracing the decision rules of the trained tree for a feature $(x_k,u_k)$ from the top node to an output leaf corresponds to the evaluation of the approximated switching conditions $\hat{g}(x_k,u_k)$.
For training, a decision tree algorithm recursively subdivides the data from the top tree node to subsequent nodes until a maximum allowable depth is reached or the number of samples at a node falls below a certain minimum. Among many other existing splitting criteria, a partitioning based on the Gini impurity measure is proposed. The decision tree algorithm starts with the initialization of the top tree node with all feature and target samples, i.e., $Q_0 = \{\xi_i,\gamma_i\}_{i=0}^{N_0-1}$ and $N_0 = N$. In general, the $m$-th node containing $N_m$ samples is referred to as $Q_m$. A split $\theta_m = (l_m, t_m)$, defined by a feature index $l_m$ and a threshold $t_m \in \mathbb{R}$, partions the node $Q_m$ in two new nodes according to $Q_m^+ (\theta_m) = \{Q_m | \xi_{i}^{(l_m)} < t_m\}$ and $Q_m^-(\theta_m) = Q_m \setminus Q_m^+(\theta_m)$, where $\xi_{i}^{(l_m)}$ denotes the $l_m$-th element of the feature vector $\xi_i = (x_i, u_i)$. The Gini impurity measure of a given node corresponds to the probability of misclassifying a randomly selected data point that is randomly labeled according to the class distribution. The impurity can be computed as
\begin{equation*}
H(Q_m) = p_{m}(\mathcal{C}_1) (1- p_{m}(\mathcal{C}_1)) + p_{m}(\mathcal{C}_2) (1- p_{m}(\mathcal{C}_2))
\end{equation*}
where  $p_{m}(\mathcal{C}_j)$ refers to the probability of selecting a data point with label $\mathcal{C}_j$ in the data set of node $Q_m$.
The impurity measure allows to determine the quality of a split $\theta_m$ with respect to a specific node $Q_m$ as
\begin{equation}
	G(Q_m,\theta_m) = \frac{N_m^+}{N_m} H(Q_m^+)  + \frac{N_m^-}{N_m}H(Q_m^-).
\end{equation}
An optimal split can therefore be determined as the impurity minimizing parameter $\theta^*_m = {\arg \min}_{\theta_m} \;G(Q_m,\theta_m)$.
As mentioned earlier, the learning algorithm recursively evaluates new optimal splits until a termination criterion such as a minimum node sample size or maximum depth is reached. The final trained decision tree provides a decision rule, starting from the top node and ending with a leaf, which can be used as an approximation to the switching condition (\ref{equ:switching_condition_ref}). {The performance of the decision tree is depicted in Fig. \ref{fig:iwp_drop_down_experiement_classified}, where the approximated switching condition $\hat{g}(x,u)$ is highlighted as black dashed line. The results of the decision tree match again almost perfectly the red reference line of the switching condition $g(x,u)$.}

\subsection{Model identification and validation}

\begin{figure}[h]
	\begin{center}
		\includegraphics[width=8.4cm]{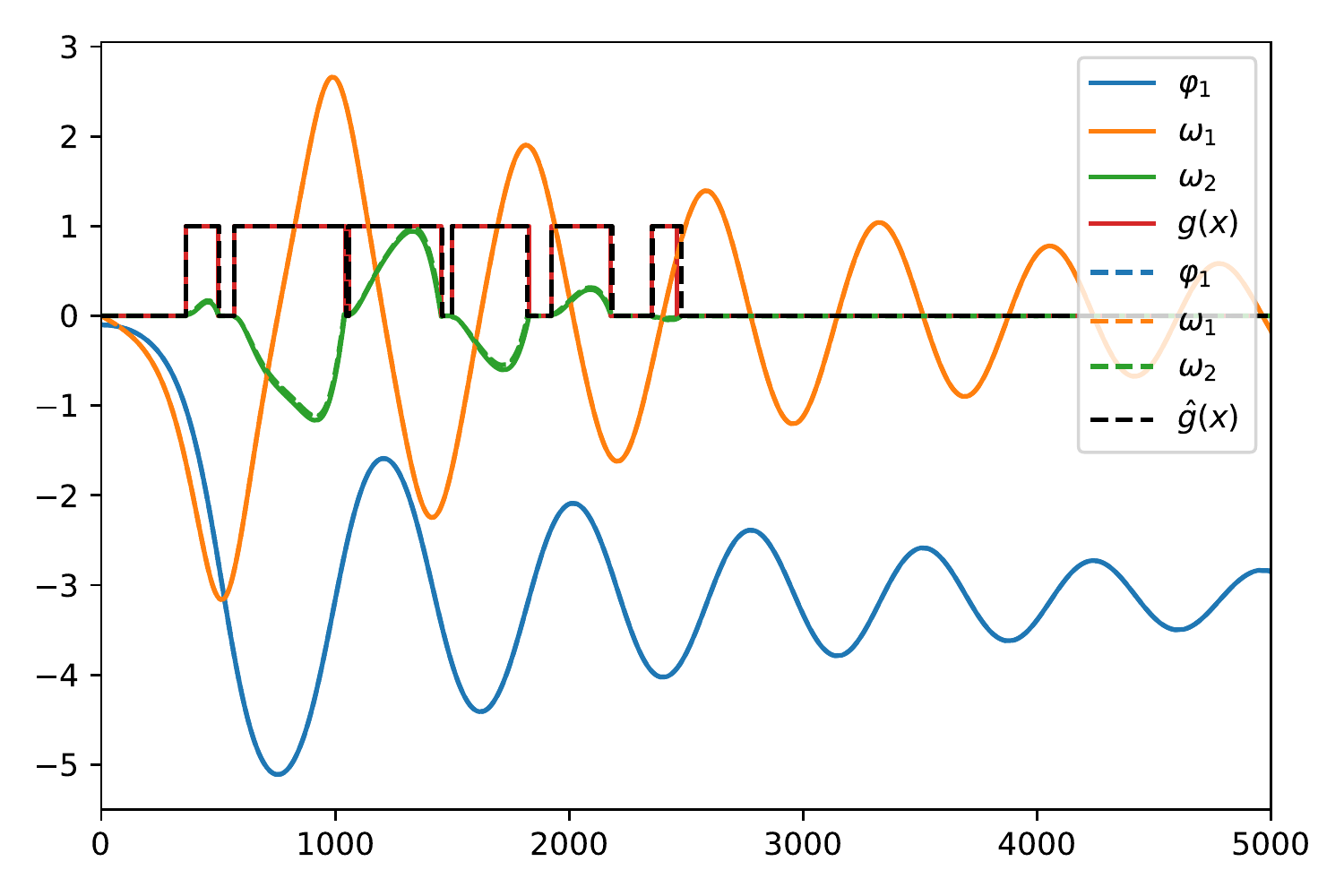}
		\caption{Validation of identified model with approximated switching condition.} 
		\label{fig:iwp_drop_down_experiement_validation}
	\end{center}
\end{figure}

The clustering step in Section 3.1 allows, in principle, an ordinary identification of two separate models belonging to $\mathcal{C}_1$ and $\mathcal{C}_2$. 
Indeed, there are several sophisticated machine learning algorithms for identifying the governing equations of nonlinear systems, see, e.g., \cite{Nelles:2001} and \cite{Brunton:2019}.  In this work, sparsity-promoting regression methods are used to identify parsimonious models. Based on the assumption that the dynamic of the recorded system is governed by a few terms (cf. the first-principle equations in Section 1), a sparse identification approach allows a kind of trade-off between model complexity and accuracy, see \cite{Brunton:2016}. For this purpose, two time-discrete models for the sticking and the non-sticking subsystem, $f_{\mathcal{C}_1}$ and  $f_{\mathcal{C}_2}$, respectively, are approximated as 
$
	x_{k+1}  =  f_{\mathcal{C}_j} (x_k, u_k) \approx \Xi_{\mathcal{C}_j} \cdot \psi(x_k,u_k),
$
where $\psi(x_k,u_k) \in \mathbb{R}^{n_\psi} $ corresponds to a vector of regression functions and $\Xi_{\mathcal{C}_j} \in \mathbb{R}^{n_x \times n_\psi}$ denotes real-valued coefficient matrices. While determining the coefficient matrices is the general objective of this step, the regressor $\psi(x_k,u_k)$ must be specified by the user through appropriate candidate functions, which should be able to correctly reflect the dynamics.
Likewise, a priori knowledge can be introduced by choosing proper regressor functions. For the following results, the candiate functions
$\psi(x_k,u_k) = [x_k, x_k^2, \sin(x_k), \cos(x_k),\mathrm{sign}(x_k),u_k]$
 were used. Note, the mathematical operations on vector $x_k$ should be interpreted as abbreviations for element-wise operations on all elements of $x_k$. As mentioned, the identification of the coefficient matrix is carried out by a sparse regression algorithm, such as, e.g., Lasso, see \cite{Brunton:2016}. Therefore, the recorded samples $x_i$, $x_{i+1}$ and $u_i$ are separated into the data matrices $X_{\mathcal{C}_j} \in \mathbb{R}^{N_{\mathcal{C}_j} \times n_x}$, $X'_{\mathcal{C}_j} \in \mathbb{R}^{N_{\mathcal{C}_j} \times n_x}$  and $U_{\mathcal{C}_j} \in \mathbb{R}^{N_{\mathcal{C}_j} \times n_u}$ with respect to their clustered labels $\gamma_i$. The estimation of the sparse coefficients $\Xi_{\mathcal{C}_j}$ is done by solving the least squares problem with $l_1$-norm penalty and weighting $\alpha \in \mathbb{R}$
{\begin{equation*}
	\Xi_{\mathcal{C}_j} = \arg \min_{\Xi_{\mathcal{C}_j}} || X_{\mathcal{C}_j}' - \Xi_{\mathcal{C}_j} \cdot \psi(X_{\mathcal{C}_j},U_{\mathcal{C}_j}) ||_2 + \alpha || \Xi_{\mathcal{C}_j} ||_1.
\end{equation*}}%
With determined coefficients, the last task of the data-driven modeling process is to combine the two identified models with the approximated switching condition, as indicated in schematic Fig. \ref{fig:inertia_wheel_pendulum}. {The validation of this approach in Fig. \ref{fig:iwp_drop_down_experiement_validation} compares the results from a drop down experiment with simulations from the identified model. The matching state trajectories show that the proposed data-driven modeling process presented in this section is indeed capable of identifying a nontrivial nonlinear system with switching structure to a high degree of accuracy. }



\section{Indirect data-driven observer design}

The data-based modeling of the inertial pendulum has been completed with the previous section.
In the present part, a practical application based on the identified model is illustrated. The objective is to design a state observer where the angle of the pendulum $\varphi_{1}$ is provided as measurement quantity. For this reason, the idea is to estimate the remaining states, including state $\omega_{2}$, which is significantly affected by the static friction, with an approximated moving horizon estimator.
Note that in \cite{Ecker:2022b} several different observers were proposed for the mathematical model of the inertial pendulum with static friction. These observers are derived from structural properties of the mathematical model, and so are not available in a pure data-driven environment without a priori knowledge. This is also a reason why numerical optimization methods such as model predictive control or moving horizon estimation are often preferred in data-driven environments.
\subsection{Moving horizon estimator}

\begin{figure}[h]
	\vspace{-0.3cm}
	\begin{center}
		\includegraphics[width=8.4cm]{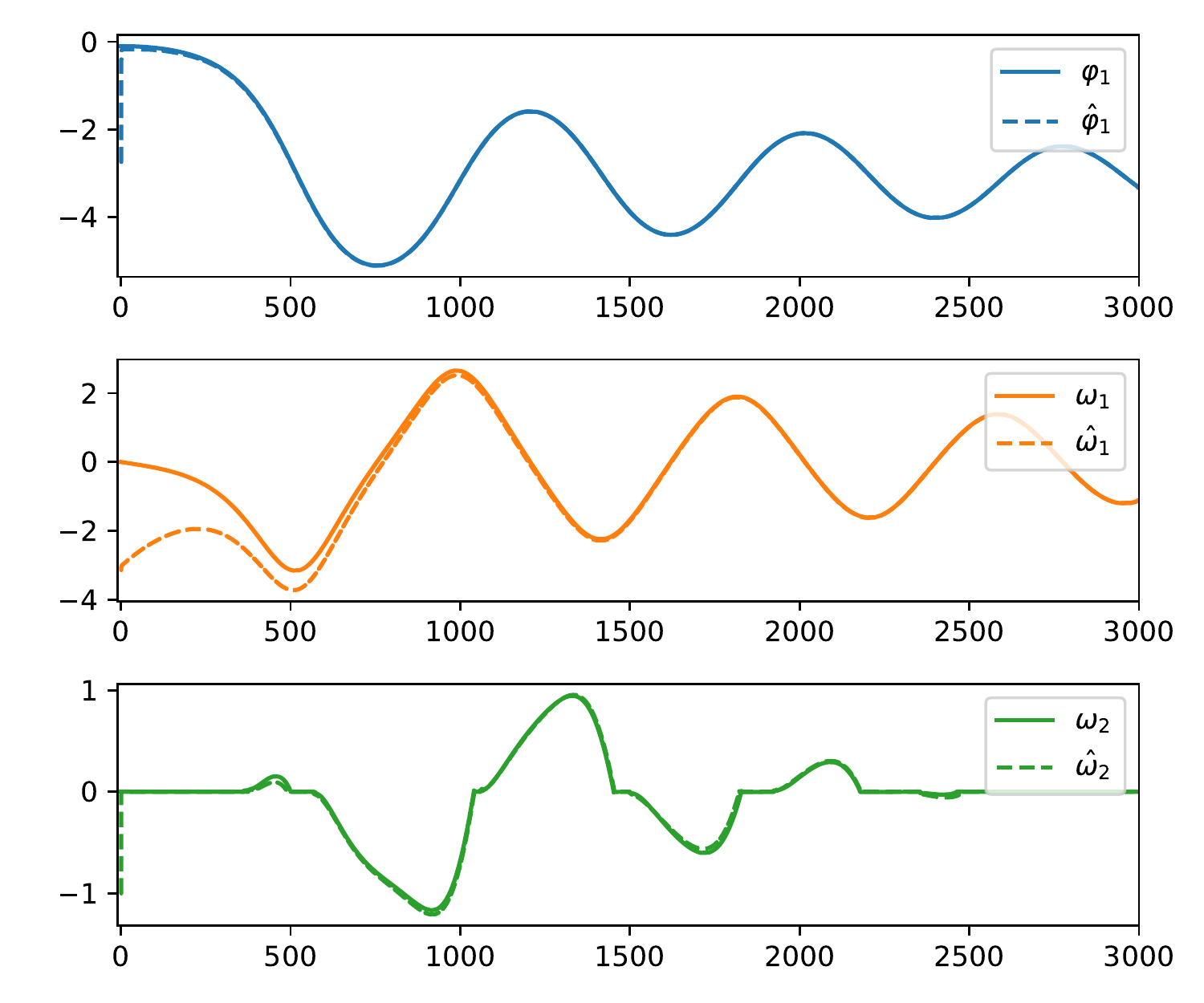}
		\caption{Validation of the moving horizon estimator.} 
		\label{fig:iwp_drop_down_experiement_mhe}
	\end{center}
\end{figure}

The principal idea of the moving horizon estimator step is to compute the most likely state values for a trajectory on a finite horizon given past measurements. In contrast to the Kalman filter, where the problem is to maximize the a posteriori estimate $p(x_{k+1}|y_k,\ldots,y_0)$, the moving horizon estimator maximizes the joint probability 
\begin{equation}
\underset{{x_{k+1}, \ldots, x_{k-N} }}{\arg \max} p(x_{k+1}, \ldots, x_{k-N} | y_{k},\ldots,y_{k-N})
\label{equ:joint_probability}
\end{equation}
on the finite horizon $N$. For the observer design problem, the previously identified models are considered with additive measurement noise $v_k$ and process disturbance $w_k$.
\begin{align*}
x_{k+1} = f_{\mathcal{C}_j}(x_k,u_k) + w_k, \qquad 
y_k = h(x_k) + v_k 
\end{align*}
The stochastic variables $v_k$ and $w_k$ are represented by their probability density functions $p(v_k)$ and $p(w_k)$.
Therefore, the maximization problem of the joint probability (\ref{equ:joint_probability}) can be written according to Bayes' rule, see \cite{Papoulis:2002}, as
{\begin{align*}
&\underset{{x_{k+1}, \ldots , x_{k-N} }}{ \max} p(x_{k+1}, \ldots, x_{k-N} | y_{k},\ldots,y_{k-N})  \\
 &= \max \prod_{i=k-N}^{k} p(y_i|x_i) p(x_{i+1}|x_i)  p(x_{k-N})
\end{align*}}%
where the conditional probability density functions $p(y_i|x_i)$ and $p(x_{i+1}|x_i)$ can be expressed in terms of the noise properties as $p(y_i|x_i) = p(v_i) \circ (y_i - h(x_i))$ and $p(x_{i+1}|x_i) = p(w_i) \circ (x_{i+1} - f(x_i,u_i))$, see, \cite{Papoulis:2002}.
Since the measurement and process noise are assumed to be zero-mean-normally distributed, i.e., $w_k \sim \mathcal{N}(0,Q_k)$ and $v_k \sim \mathcal{N}(0,R_k)$ with covariance matrices $Q_k \in \mathbb{R}^{3 \times 3} > 0$ and $R_k \in \mathbb{R} > 0$, it is convenient to maximize the logarithm of the joint probability density, which is valid due to the monotonicity of the transformation.
{
\begin{align*}
 &\max \sum_{i=k-N}^{k} \log p(y_i|x_i) + \log p(x_{i+1}|x_i) + \log p(x_{k-N})\\
 = &\min \sum_{i=k-N}^{k} ||\underbrace{y_i -h(x_i)}_{v_i}||^2_{R_i^{-1}} + ||\underbrace{x_{i+1} -f(x_i,u_i))}_{w_i}||^2_{Q_i^{-1}}  \\ &\qquad \qquad \qquad \qquad \qquad \qquad \qquad \quad + ||x_{k-N} -\hat{x}_{k-N}||^2_{P_i^{-1}} 
\end{align*}}%
The logarithm applied to the normally distributed probability densities results in a minimization problem due to a sign change. Furthermore, the variance terms are abbreviated according to the $P$-weighted squared norm $||x||_{P}^2 = x^T P x$. The optimization problem of the moving horizon estimator can therefore be be formulated as
{
\begin{align*}
\hat{x}_{k-N},\ldots,\hat{x}_{k+1} &=\underset{{x_{k-N},\ldots,x_{k+1}}}{\arg \min} || x_{k-N} - \hat{x}_{k-N}||^2_{P_i^{-1}}  \\ & \qquad +\sum_{i=k-N}^{k} \left(||v_i||^2_{R_{i}^{-1}} + ||w_i||^2_{Q_i^{-1}}\right) & \\
\mathrm{s.t. }  \quad  x_{k+1} & =  \left\{ \begin{array}{ll}
\hat{f}_{C_1} (x_k, u_k) + w_k & C_1 =  \hat{g}(x_k, u_k) \\
\hat{f}_{C_2}(x_k, u_k) + w_k & C_2 =  \hat{g}(x_k, u_k)
\end{array} \right. \\
y_k  & = \varphi_{1,k} + v_k 
\end{align*}
 }%
{A special aspect of the moving horizon estimator is the combination of the general optimization idea with the switching structure of the identified model in the constraints. }
The estimation performance of the proposed indirect data-driven observer with a finite horizon $N=10$ and sampling time $T_a = 5\, \text{ms}$ is presented in Fig. \ref{fig:iwp_drop_down_experiement_mhe}, where recorded state trajectories (solid lines) are compared to the moving horizon estimates (dashed lines) for a pendulum drop down experiment with an initial observer error. The results show that the initial observer error converges to zero and that the data-driven observer is able to handle both the switching dynamics and the effects of static friction on the inertia wheel. 
\section{Conclusion}
In this contribution, an indirect-data driven observer design for the inertia wheel pendulum with static friction has been presented. The data-driven modeling of the system under consideration with inherent structural change is composed of a mixture of unsupervised and supervised machine learning algorithms that allow the two dynamics of the system to be clustered and classified. Combining the two identified models with an approximated switching condition yields the required model accuracy. The subsequent observer is based on the moving horizon prediction. Further research focuses on the improvement of algorithms for highly nonlinear systems and the automated cluster order determination while taking model accuracy into account.

\bibliography{ifacconf}             

\end{document}

%% file: figures/schematic_and_image_inertia_wheel_pendulum_en-min.pdf_tex
\begingroup%
  \makeatletter%
  \providecommand\color[2][]{%
    \errmessage{(Inkscape) Color is used for the text in Inkscape, but the package 'color.sty' is not loaded}%
    \renewcommand\color[2][]{}%
  }%
  \providecommand\transparent[1]{%
    \errmessage{(Inkscape) Transparency is used (non-zero) for the text in Inkscape, but the package 'transparent.sty' is not loaded}%
    \renewcommand\transparent[1]{}%
  }%
  \providecommand\rotatebox[2]{#2}%
  \newcommand*\fsize{\dimexpr\f@size pt\relax}%
  \newcommand*\lineheight[1]{\fontsize{\fsize}{#1\fsize}\selectfont}%
  \ifx\svgwidth\undefined%
    \setlength{\unitlength}{255.11811024bp}%
    \ifx\svgscale\undefined%
      \relax%
    \else%
      \setlength{\unitlength}{\unitlength * \real{\svgscale}}%
    \fi%
  \else%
    \setlength{\unitlength}{\svgwidth}%
  \fi%
  \global\let\svgwidth\undefined%
  \global\let\svgscale\undefined%
  \makeatother%
  \begin{picture}(1,0.77777778)%
    \lineheight{1}%
    \setlength\tabcolsep{0pt}%
    \put(-6.04759598,2.23992642){\color[rgb]{0,0,0}\makebox(0,0)[lt]{\lineheight{0}\smash{\begin{tabular}[t]{l} \end{tabular}}}}%
    \put(0,0){\includegraphics[width=\unitlength,page=1]{figures/schematic_and_image_inertia_wheel_pendulum_en-min.pdf}}%
    \put(0.5252379,0.49314512){\color[rgb]{0,0,0}\makebox(0,0)[lt]{\lineheight{283.8500061}\smash{\begin{tabular}[t]{l}$\mathrm{Pendulum}$\end{tabular}}}}%
    \put(0,0){\includegraphics[width=\unitlength,page=2]{figures/schematic_and_image_inertia_wheel_pendulum_en-min.pdf}}%
    \put(0.52629714,0.54453222){\color[rgb]{0,0,0}\makebox(0,0)[lt]{\lineheight{283.8500061}\smash{\begin{tabular}[t]{l}$\mathrm{Motor}$\end{tabular}}}}%
    \put(0.52629714,0.59591941){\color[rgb]{0,0,0}\makebox(0,0)[lt]{\lineheight{283.8500061}\smash{\begin{tabular}[t]{l}$\mathrm{Wheel}$\end{tabular}}}}%
    \put(0.71457197,0.04106451){\color[rgb]{0,0,0}\makebox(0,0)[lt]{\lineheight{283.8500061}\smash{\begin{tabular}[t]{l}$\mathrm{Counterbalance}$\end{tabular}}}}%
    \put(0,0){\includegraphics[width=\unitlength,page=3]{figures/schematic_and_image_inertia_wheel_pendulum_en-min.pdf}}%
  \end{picture}%
\endgroup%

%% file: figures/iwp_switched_system.pdf_tex
\begingroup%
  \makeatletter%
  \providecommand\color[2][]{%
    \errmessage{(Inkscape) Color is used for the text in Inkscape, but the package 'color.sty' is not loaded}%
    \renewcommand\color[2][]{}%
  }%
  \providecommand\transparent[1]{%
    \errmessage{(Inkscape) Transparency is used (non-zero) for the text in Inkscape, but the package 'transparent.sty' is not loaded}%
    \renewcommand\transparent[1]{}%
  }%
  \providecommand\rotatebox[2]{#2}%
  \newcommand*\fsize{\dimexpr\f@size pt\relax}%
  \newcommand*\lineheight[1]{\fontsize{\fsize}{#1\fsize}\selectfont}%
  \ifx\svgwidth\undefined%
    \setlength{\unitlength}{226.77165354bp}%
    \ifx\svgscale\undefined%
      \relax%
    \else%
      \setlength{\unitlength}{\unitlength * \real{\svgscale}}%
    \fi%
  \else%
    \setlength{\unitlength}{\svgwidth}%
  \fi%
  \global\let\svgwidth\undefined%
  \global\let\svgscale\undefined%
  \makeatother%
  \begin{picture}(1,0.625)%
    \lineheight{1}%
    \setlength\tabcolsep{0pt}%
    \put(0,0){\includegraphics[width=\unitlength,page=1]{iwp_switched_system.pdf}}%
    \put(0.49012527,0.48446903){\makebox(0,0)[lt]{\smash{\begin{tabular}[t]{l}$\hat{g}(x_k,u_k)$\end{tabular}}}}%
    \put(0.48597207,0.32612009){\makebox(0,0)[lt]{\smash{\begin{tabular}[t]{l}$\hat{f}_{\mathcal{C}_1}(x_k,u_k)$\end{tabular}}}}%
    \put(0.48124294,0.15762057){\makebox(0,0)[lt]{\smash{\begin{tabular}[t]{l}$\hat{f}_{\mathcal{C}_2}(x_k,u_k)$\end{tabular}}}}%
    \put(0,0){\includegraphics[width=\unitlength,page=2]{iwp_switched_system.pdf}}%
    \put(0.00008949,0.3187757){\makebox(0,0)[lt]{\smash{\begin{tabular}[t]{l}$x_k,u_k$\end{tabular}}}}%
  \end{picture}%
\endgroup%